\documentclass[12pt,reqno]{amsart}

\usepackage{amsmath,amsfonts,amsbsy,amsgen,amscd,mathrsfs,amssymb,amsthm}
\usepackage{enumerate,mathtools}

\usepackage[usenames,dvipsnames]{xcolor}
\usepackage[colorlinks=true,citecolor=blue,linkcolor=blue]{hyperref}

\usepackage{tikz}
\usetikzlibrary{shadings}

\newtheorem{thm}{Theorem}[section]

\newtheorem{prop}[thm]{Proposition}

\newtheorem{rem}[thm]{Remark}

\newcommand{\pd}{\partial}

\theoremstyle{definition}

\numberwithin{equation}{section} 
\numberwithin{figure}{section}
\numberwithin{table}{section}

\newcommand{\bE}{\mathbf{E}}
\newcommand{\eps}{\varepsilon}
\newcommand{\bR}{\mathbf{R}}

\newcommand{\bP}{\mathbf{P}}
\newcommand{\Om}{\Omega}

\newcommand{\la}{\lambda}

\newcommand{\vf}{\varphi}

\begin{document}

\title[Vector-valued functions measure concentration ]{Measure concentration for vector valued functions on Hamming cube}



\author{Alexander Borichev}
\address{(A.B.) Department of Mathematics, Aix-Marseille University, CNRS, Centrale Marseille, I2M, France}
\email{alexander.borichev@math.cnrs.fr}

\author{Alexander Volberg}
\address{(A.V.) Department of Mathematics, MSU, 
East Lansing, MI 48823, USA }
\email{volberg@math.msu.edu}

\begin{abstract}
We prove here the concentration of measure inequality for ``Lipschitz'' function on the Hamming cube
with values in any Banach spaces of finite cotype.
 \end{abstract}

\thanks{ The research of the  authors is supported  by 
Hausdorff Research Institute for Mathematics in Bonn, during the trimester
program ``Boolean Analysis in Computer Science", funded by the Deutsche
Forschungsgemeinschaft (DFG, German Research Foundation) under Germany's
Excellence Strategy – EXC-2047/1 – 390685813; AV NSF grant DMS-2154402 }



\maketitle


\section{Introduction}
\label{intro}
Measure concentration results and related sharp estimates of moments of function of (independent) random variables 
enjoyed a constant attention of many mathematicians, the bibliography is immense. We refer here only to some texts that contain lots of further references: \cite{HVNVW2}, \cite{L}, \cite{LO}, \cite{AW}, \cite{T1}.

\medskip

Hamming cube is $\Om_2^n:=\{-1, 1\}^n$. Integration with respect to uniform probability measure is denoted by $\bE$.  Any function
on $\Om_2^n$ is a multilinear polynomial. Symbol $\pd_i$ denotes the usual partial derivative, $i=1, \dots, n$. Also
$$
D_i f= x_i\pd_i f(x)\,.
$$
This operation can be written as
$$
D_i f(x) =\frac{f(x)- f(x_1, \dots, x_{i-1}, -x_i, x_{i+1}, \dots, x_n)}{2}\,.
$$
All $f$ below are
$$
f :\Om_2^n \to E,
$$
where $E$ is a normed space.

We are proving here some estimates for functions on Hamming cube with values in any Banach space of finite cotype that imply measure concentration results for all spaces of finite cotype.

\section{Main ideas}
\label{idea}

The main ideas are to use Gross log-Sobolev inequality (LSI) combining with a beautiful observation from Cordero-Erasquin--Eskenazis \cite{CE}, but
applied to a different object. Then we compare two differential inequalities not unlike this is done in the famous Herbst argument. The use of Proposition 4.2 of \cite{IVHV} finishes the proof of Theorem \ref{Ke} below.
\begin{thm}
\label{Ke}
For any Banach space $E$ with cotype $Q$ and constant of cotype $C(E)$ the following measure concentration result holds under the Lipschitz assumption \eqref{lip2}:
$$
\bE e^{\frac1{4e} \|f-\bE f\|_E^2} \le e^{c_0 C(E)^2 Q^2},
$$
where $c_0$ is a universal constant.
\end{thm}
\begin{rem} 
After proving Theorem \ref{Ke} we found out that it has been also observed by Ramon van Handel.
\end{rem}
\begin{rem} 
A feature of the theorem that can be emphasized is the independence of the constant in the exponent of either the cotype or cotype constant.
\end{rem}

\section{Combining  LSI \cite{G} and Cordero-Erasquin--Eskenazis \cite{CE}}
\label{comb}

Let $E$ be a normed space. Let $\psi_2(t):= t^2 \log (e+t^2)$ be the function giving us the Orlicz space $L^2 \log L$, where the norm of a test function $g\ge 0$ is given by
$$
\|g\|_{L^2 \log L}:= \inf\{\la: \bE \psi_2\big( \frac{g}{\la}\Big) \le 1\}\,.
$$

\begin{prop}(Inequality (11) of \cite{CE})
\label{ent}
Let $g$ be a nonnegative function on $\Om_2^n$. Then
$$
\bE g^2(\log g^2 -\log \bE g^2) \le 2 \|g\|_{L^2\log L}^2\,.
$$
\end{prop}

\medskip

For function $g: \Om_2^n\to \bR$ let 
$$
Mg (x)  :=\Big(\sum_{i=1}^n [\big( D_i g(x)\big)_+]^2\Big)^{1/2}\,.
$$
The following result was proved  in \cite{T} (Theorem 1.4).

\begin{prop}(Inequality (19) of \cite{CE})
\label{M}
Let $g$ be a nonnegative function on $\Om_2^n$. Let $\bP\{ x: g(x)=0\}\ge \frac12$. Let $1<p<\infty$. Then
$$
 \|g\|_{L^p\log^{p/2} L}^p \le \kappa_p\bE\, (Mg)^p\,.
 $$
\end{prop}

\begin{rem}
The sharp $\kappa_p$ is not known. We will use this only for $p=2$, see Proposition \ref{entGT}. It follows immediately from Proposition \ref{entG}.
\end{rem}

\medskip

\begin{prop}(Log-Sobolev inequality on Hamming cube)
\label{entG}
Let $g$ be a nonnegative function on $\Om_2^n$. Then
$$
\bE g^2(\log g^2 -\log \bE g^2) \le 2 \|\nabla g\|_{L^2}^2\,.
$$
\end{prop}

\begin{prop}
\label{entGT}
Let $g$ be a nonnegative function on $\Om_2^n$. Then
$$
\bE g^2(\log g^2 -\log \bE g^2) \le 4 \bE (Mg)^2\,.
$$
\end{prop}

Actually it is only this proposition that we will be using below. But it might happen that using propositions due to Talagrand one can maybe improve some constants below.

\bigskip

Let now $f(x)$ be a function on the Hamming cube with values in a normed space $E$, and let $F(x)$ be its extension
by the same multilinear formula to the whole $\bR^n$.

Let $p$ be a large number and $g(x) := \|f(x)\|_E^{p/2}$. We wish to see what is $Mg$, and, hence, what is  $(D_i g)_+(x), x\in \Om_2^n$.
Following \cite{CE} let $e_i=(0...0,1,0...0)$ and
$$
\vf_{x, i} (s) := F(x+ s\,e_i),\quad x\in \Om_2^n, \,\, s\in \bR,
$$
it is a linear function, therefore,
$$
\frac{\vf_{x, i} (s_1)  + \vf_{x, i} (s_2) }{2} = \vf_{x, i}\Big(\frac{s_1+s_2}{2}\Big)\,.
$$
Hence  $s\to \|\vf_{x,i}(s)\|_E$ is convex, and as $t\to t^{p/2}$ is increasing, we get
$$
\|\vf_{x, i}\Big(\frac{s_1+s_2}{2}\Big)\|_E^{p/2} \le \Big(\frac{\|\vf_{x, i}(s_1)\|_E + \|\vf_{x, i}(s_2)\|_E}{2}\Big)^{p/2}
$$
As $t\to t^{p/2}$ is  convex $p\ge 2$ we finally get
\begin{equation}
\label{conv}
\|\vf_{x, i}\Big(\frac{s_1+s_2)}{2}\Big)\|_E^{p/2} \le \frac{\|\vf_{x, i}(s_1)\|_E^{p/2} + \|\vf_{x, i}(s_2)\|^{p/2}_E}{2},
\end{equation}
which is the separate  convexity of $\|F(x)\|_E^{p/2}$ in every variable in $\bR^n$.

\bigskip

This allows us to estimate $(D_i \|f(x)\|_E^{p/2})_+(x), x\in \Om_2^n$, as it has been done in a similar situation in \cite{CE}.
In fact, assume first that $x_i=1$. If $\|f(....,1,...)\|_E \le \|f(....,-1,...)\|_E$ then $(D_i \|f(x)\|_E^{p/2})_+(x)=0$. 

If $\|f(....,1,...)\|_E > \|f(....,-1,...)\|_E$ then 
$$
(D_i \|f(x)\|_E^{p/2})_+(x) = \frac{\|F(....,1,...)\|_E^{p/2}-\|F(....-1,....)\|_E^{p/2}}{2},
$$ 
which the slope of the chord of the graph of convex function on the interval $[-1,1]$,
and the value at the left point is smaller than the value at the right  point of the interval over which the chord lies. But then convexity says that the slope
over any interval $[1, 1+\eps], \eps>0$, will be at least as big.
 We conclude that
 in this case
 \begin{equation}
 \label{deriv1}
 (D_i \|f(x)\|_E^{p/2})_+(x) \le \Big|\frac{\partial \big(\|F(x)\|_E^{p/2}\big)}{\partial x_i} (x)\Big|, \quad x=(.....,1,....)\,.
 \end{equation}
By what we said before the same is automatically true for the case $\|f(....,1,...)\|_E \le \|f(....,-1,...)\|_E$, as  then the left hand side is zero.

\bigskip

Now assume  that $x_i=-1$. If $\|f(....,-1,...)\|_E \le \|f(....,1,...)\|_E$ then $(D_i \|f(x)\|_E^{p/2})_+(x)=0$.

If $\|f(....,-1,...)\|_E > \|f(....,1,...)\|_E$ then 
$$
(D_i \|f(x)\|_E^{p/2})_+(x) = \frac{\|F(....,-1,...)\|_E^{p/2}-\|F(....1,....)\|_E^{p/2}}{2},
$$ 
which the absolute value of the slope of the chord of the graph of convex function on  interval $[-1,1]$, 
but this time the value of our convex function on the left end-point is bigger than on the right end-point.

Then convexity says that the absolute value of the slope
over any interval $[-1-\eps, -1], \eps>0$, will be at least as big.
 We conclude that
 in this case
 \begin{equation}
 \label{deriv2}
 (D_i \|f(x)\|_E^{p/2})_+(x) \le \Big|\frac{\partial \big(\|F(x)\|_E^{p/2}\big)}{\partial x_i} (x)\Big|, \quad x=(.....,-1,....)\,.
 \end{equation}
 
 Finally, we see that we always have
 \begin{equation}
 \label{deriv}
\!\!\!\!\! (D_i \|f(x)\|_E^{p/2})_+(x) \le \Big|\frac{\partial \big(\|F(x)\|_E^{p/2}\big)}{\partial x_i} (x)\Big|, \quad x\in \Om_2^n, i=1, \dots, n\,.
 \end{equation}

\bigskip

Now let us apply  Proposition \ref{entGT} and \eqref{deriv}. Then for 
$$
g=\|f(x)\|_E^{p/2}
$$
 we get \eqref{Orlicz0} below, 
\begin{equation}
\label{Orlicz0}
\bE g^2(\log g^2 -\log \bE g^2) \le 4 \bE (Mg)^2
\end{equation}

So, 
\begin{align}
\label{Orlicz1}
& \bE \|f\|_E^{p}(\log \|f\|_E^{p} -\log \bE \|f\|_E^{p}) \le              \notag
\\
&4\Big[ \bE \Big( \sum_{i=1}^n  \Big|\frac{\partial \big(\|F(x)\|_E^{p/2}\big)}{\partial x_i} (x)\Big|^2\Big)\Big]^{1/2} \,.
\end{align}

\medskip


\bigskip

Now we can apply inequality (34)  of \cite{CE} under the assumption of Lipschitzness:
\begin{align}
\label{Lip}
&\bE \Big( \sum_{i=1}^n  \Big|\frac{\partial \big(\|F(x)\|_E^{p/2}\big)}{\partial x_i} (x)\Big|^2\Big) \le          \notag
\\
&\frac{p^2}{4} \bE_x \Big(\|f(x)\|_E^{p-2} \bE_\delta \Big\|\sum_{i=1}^n \delta_i D_i f(x)\Big\|^2_E\Big) \le \frac{p^2}{4} \Big[\bE_x \big(\|f(x)\|_E^{p}\big)\Big]^{1-\frac2{p}} \,.
\end{align}
The last inequality is valid if we assume the {\it Lipschitz} property of $f$:
\begin{equation}
\label{lip1}
P^2(f) :=\bE_\delta \Big\|\sum_{i=1}^n \delta_i D_i f(x)\Big\|^2_E \le 1, \quad \forall x\in \Om_2^n\,.
\end{equation}

\medskip 

If we do not assume \eqref{lip1} we still have
\begin{align}
\label{Grad}
&\bE \Big( \sum_{i=1}^n  \Big|\frac{\partial \big(\|F(x)\|_E^{p/2}\big)}{\partial x_i} (x)\Big|^2\Big) \le \frac{p^2}{4} \bE_x \Big(\|f(x)\|_E^{p-2}\, \bE_\delta \Big\|\sum_{i=1}^n \delta_i D_i f(x)\Big\|^2_E\Big) \le              \notag
\\
& \frac{p^2}{4} \Big[\bE_x \big(\|f(x)\|_E^{p}\big)\Big]^{1-\frac2{p}}  \cdot  \Big\{\bE_x\Big(\bE_\delta \Big\|\sum_{i=1}^n \delta_i D_i f(x)\Big\|^2_E\Big)^{p/2}\Big\}^{2/p} \,.
\end{align}

\medskip

The natural Lipschitz norm in the vector-valued case would be the one in the left hand side of \eqref{lip0} below. It is called {\it weak Lipschitz norm}.
Notice that weak Lipschitz property in \eqref{lip0} will also give  estimate  \eqref{Lip} (this time by inequality (32) of \cite{CE}).
\begin{equation}
\label{lip0}
\sup_{x} \sup_{\|\xi\|_{E^*}=1} \sum_{i=1}^n \langle \xi, D_if(x)\rangle^2 \le 1,
\end{equation}
where $\xi$  is in the sphere of the dual $E^*$. This is a smaller gradient than in \eqref{lip1} because
\begin{align*}
&\sup_{\|\xi\|_{E^*}=1} \sum_{i=1}^n \langle \xi, D_if(x)\rangle^2= \sup_{|\xi\|_{E^*}=1} \bE_\delta\Big|\sum\delta_i  \langle \xi, D_if(x)\rangle\Big|^2 \le 
\\
&\bE_\delta\Big|\sup_{|\xi\|_{E^*}=1}\sum\delta_i  \langle \xi, D_if(x)\rangle\Big|^2= \bE_\delta\Big\| \sum_{i=1}^n \delta_i  D_if \Big\|_E^2\,.
\end{align*}

\begin{rem}
\label{gradi}
The fact that the gradient in \eqref{lip0} is the smallest can be used in estimate \eqref{Orlicz2} below.
 But looks like that for the final result  of Theorem \ref{KE} we still need a bigger gradient from \eqref{lip1}. 
 This is because Proposition 4.2 of \cite{IVHV} (Pisier--Poincar\'e inequality in spaces of finite cotype) uses gradient from \eqref{lip1}. I am grateful for this remark to Alexandros Eskenazis. 
\end{rem}

\medskip

Finally, there is the yet another gradient, which can be used to define vector-valued Lipschitz functions:

\begin{equation}
\label{lip2}
\Gamma^2(f):=\sum_{i=1}^n \Big\|D_i f(x)\Big\|^2_E\ \le 1, \quad \forall x\in \Om_2^n\,.
\end{equation}

\medskip

Considering all that, we finally get that under either \eqref{lip0} or \eqref{lip1} or \eqref{lip2} 
(and the gradient in \eqref{lip0} is the smallest, and, so, the best for our purpose) the following holds
\begin{equation}
\label{Orlicz2}
\bE \|f\|_E^p (\log \|f\|_E^p -\log \bE \|f\|_E^{p}) \le p^2 \Big(\bE_x \big(\|f(x)\|_E^{p}\big)\Big)^{1-\frac2{p}} \,.
\end{equation}





We denote
$$
a(p):= \bE\|f\|_E^p\,.
$$
So,
$a'(p) = \bE\big(\|f\|_E^p \log\|f\|_E\big)$.

\medskip

Then \eqref{Orlicz2} reads as follows
\begin{equation}
\label{diff1}
a'(p)  \le \frac1p a(p)\log a(p) +  p \,a(p)^{1-\frac2{p}} \,.
\end{equation}

Consider now
$$
\beta(p) :=\frac{\log a(p)}p\,.
$$
By \eqref{diff1} we have
\begin{multline*}
0\le \beta'(p)=\frac{a'(p)}{pa(p)}-\frac{\log a(p)}{p^2}
\\ \le 
\frac{\log a(p)}{p^2}+\frac{1}{a(p)^{2/p}}
-\frac{\log a(p)}{p^2}
\\=
e^{-2\beta(p)}
\end{multline*}

On the other hand, \cite{IVHV} tells us that if we are in the Lipschitz situation
\eqref{lip1}:
\begin{equation}
\label{lip2}
P^2(f) :=\bE_\delta \Big\|\sum_{i=1}^n \delta_i D_i f(x)\Big\|^2_E \le 1, \quad \forall x\in \Om_2^n\,.
\end{equation}
then
\begin{equation}
\label{ivhv}
\beta(p)  \le \log C(E) +\log p + \frac32 \log_+ \frac{Q}{p}, \quad p\in [1, \infty),
\end{equation}
where $C(E)$ is a constant of cotype of $E$, and $Q$ is the cotype of $E$.

Function $\beta$ is increasing, but the right hand side function is decreasing in $p\in [1, Q]$. Henceforth, in the Lipschitz situation \eqref{lip2} \eqref{ivhv} immediately self-improves to
\begin{equation}
\label{ivhv1}
\beta(p)  \le \log C(E) +\log Q, \quad p\in [1, Q]\,.
\end{equation}

\subsection{Estimate of $p$-th moments for general Banach space of finite cotype $Q$}
\label{Egeneral}

In general we have
\begin{equation}
\label{ivhv2}
e^{\beta(p)}  \le C(E) Q, \quad p\in [1, Q], 
\end{equation}
and
\begin{equation}
\label{bega}
\beta(x) \le \gamma(x) := \frac12 \log (2x +C(E)^2 Q^2 - 2Q),\quad x\ge Q\,. 
\end{equation}

Indeed, $\gamma' = e^{-2\gamma}$, and $\beta(Q)\le \gamma(Q)$. If for some $Q<x<y$ we have 
$$
\beta(x)=\gamma(x)
$$
 and  for $t: x<t<y$ $\beta(t)>\gamma(t)$, then
$\beta'(t) \le e^{-2\beta(t)}  < e^{-2 \gamma(t)} =\gamma'(t)$.
Then 
$$
\beta(t) < \gamma(t), t\in (x, y),
$$
which contradicts our assumption  above.  So \eqref{bega} is proved.

\begin{thm}
\label{sqrtp}
Let $E$ be a Banach space of finite cotype $Q$ and constant $C(E)$.  Let Lipschitz condition \eqref{lip2} hold. Then
$$
(\bE \|f-\bE f\|_E^p)^{1/p} \le  \begin{cases} 
& C\sqrt{2p+C(E)^2 Q^2 -2Q}, \quad p\ge Q,
\\
& C(E)Q, \quad 1\le p\le Q\,.
\end{cases}
$$
In particular, for large $p$ the growth is of the order $\sqrt{p}$ for any Banach space of finite cotype.
\end{thm}

\bigskip

To estimate 
$$
\bE e^{\tau \|f-\bE f\|_E^2}
$$ 
with small universal $\tau$, we need to estimate
\begin{equation}
\label{sigma}
\Sigma:= \sum_{n=0}^\infty \frac{\tau^n a(2n)}{n!} \asymp \sum_n e^{n \log \tau + 2n \beta(2n)-n\log n +n -\frac12 \log n}\,.
\end{equation}
Hence,
$$
\Sigma \asymp 1+ \sum_{n=1}^{Q/2}  e^{n\log\tau + 2n \log (C(E)Q)-n\log n +n } + 
$$
$$
 \sum_{n=Q/2+1}^{\infty}  e^{n\log\tau + n \log (4n+C(E)^2 Q^2 -2Q) -n\log n +n }\,.
$$
Function $n\to n\log\tau + 2n \log (C(E)Q)-n\log n +n$ first grows, then decreases, but this happens only after $n\approx C(E)^2 Q^2 \ge Q^2$, 
so we can estimate all terms of the first sum by the term at $n=Q/2$.

So, if we choose 
$$
\tau\le 1/e
$$
we get
$$
\Sigma_1 \le Q e^{Q\log (2C(E))}\,.
$$
Also if we choose
$$
\tau\le 1/4e,
$$
then
$$
\Sigma_2 =\sum_{n=Q/2+ 1}^{C(E)^2 Q^2}\dots + \sum_{n=C(E)^2 Q^2+1}^{\infty}\dots\le
$$
$$
5C(E)^2 Q^2 e^{C(E)^2 Q^2} + C_0 e^{C(E)^2 Q^2},
$$
with universal $C_0$ (we just estimate the second sum using a geometric progression).

Therefore we proved
\begin{thm}
\label{KE}
For any Banach space $E$ with cotype $Q$ and constant of cotype $C(E)$ the following measure concentration result holds under the Lipschitz assumption \eqref{lip2}:
$$
\bE e^{\frac1{4e} \|f-\bE f\|_E^2} \le e^{c_0 C(E)^2 Q^2},
$$
where $c_0$ is a universal constant.
\end{thm}

\subsection{Let us discuss the sharpness}
Notice that $\Sigma$ from \eqref{sigma} has the following estimate from below (say, $\tau$ is not very small here):
$$
\!\!\!\!\log \Sigma  \ge \tau Q^2 \log \tau + 2\tau Q^2\log Q - \tau Q^2 \log (\tau Q^2) + \tau Q^2 - \log \frac{Q}{2} \ge \frac{\tau}{2}\, Q^2.
$$
We just took one term for $n= \tau Q^2$ in the sum $\Sigma$ in \eqref{sigma}, and we used the assumption $\beta (\tau Q^2) = \log Q$. This assumption means

\begin{equation}
\label{logQ}
(\bE\|f\|^p_E)^{1/p} =Q, \quad p=\tau Q^2\,.
\end{equation}

This should hold for a certain $f:\Om_2^n \to E$ such that $\bE f=0$ and $\big\| \bE_\delta \|\sum_{i=1}^n \delta_i D_i f(x)\|_E^2\big\|_{L^\infty (\{-1,1\}^n)} \le 1$, which is our Lipschitz assumption \eqref{lip1}.

So the question is whether we can reconcile the assumption \eqref{logQ} with this Lipschitz assumption by presenting a corresponding $f$.






\section{An application. A comparison with \cite{HT}}
\label{matrix1}

There are two different Lipschitz conditions for spaces 
$$
E=(M_{d\times d}, \|\cdot\|_{S_p}),
$$
where $Q=p$, $C(E)=1$, \cite{HVNVW2}, Proposition 7.1.11. The first type of Lipschitz condition is \eqref{lip2}, that is
\begin{equation}
\label{lip3}
P^2_p(f) :=\Big(\bE_\delta \Big\|\sum_{i=1}^n \delta_i D_i f(x)\Big\|^2_{S_p}\Big)^{1/2} \le 1, \quad \forall x\in \{-1,1\}^n.
\end{equation}
Another type of Lipschitz condition is $\forall x\in \{-1,1\}^n$
\begin{equation}
\label{LipM}
\!\!\!\!\!\!K^2_p :=\Big\|\Big(\sum_{i=1}^n  D_i f(D_i f)^*(x)\Big)^{1/2}_{S_p}\Big\|+ \Big\|\Big(\sum_{i=1}^n  (D_i f)^*D_i f(x) \Big)^{1/2}\Big\|_{S_p} \le 1.
\end{equation}
For $p=\infty$ the norm is the operator norm as usual.

\bigskip

The relationships between $P^2_p$ and $K^2_p$ is via non-commutative Khintchine inequality, namely (below $c_1, c_2$ are absolute constants).
\begin{equation}
\label{Kh}
c_1 \, K^2_p \le P^2_p \le c_2\, \min[\sqrt{p}, \sqrt{\log d}]\, K^2_p, \quad  2\le p \le \infty\,.
\end{equation}
See \cite{P}, page 106.

\bigskip

Using \eqref{ivhv2}, let us consider the example of $E=(M_{d\times d}, \|\cdot\|_{op})$ , where $Q=\log d$ and $C(E)=1$.
Then the previous display establishes that if 
$$
\bE_\delta \Big\|\sum_{i=1}^n \delta_i D_i f(x)\Big\|^2_{op} \le 1, \quad \forall x\in \Om_2^n,
$$
then (with $\bE f=0$)
\begin{equation}
\label{m1} 
\big(\bE\|f\|_{op}^p\big)^{1/p} = a(p)^{1/p} =e^{\beta(p)}  \le^{\eqref{ivhv1}} \log d, \quad p\in [1, \log d]\,.
\end{equation}
This is worse than Theorem 3.2 of \cite{HT}, where the right hand side is of the order $\sqrt{\log d}$ for $p=1$.

\medskip

If we consider $E=(M_{d\times d}, \|\cdot\|_{S_p})$, $2\le p\le \log d$, where $S_p$ is a Schatten--von Neumann class, $Q=p$ and $C(E)=1$,
we get
\begin{thm}
Let $f:\{-1,1\}^n\to M_{d\times d}$, and $E f=0$. Then
\label{P2p}
\begin{align}
\label{m2} 
&\!\!\!\!\!P^2_p\le 1 \Rightarrow \big(\bE\|f\|_{S_p}^p\big)^{1/p} = a(p)^{1/p}=e^{\beta(p)}  \le^{\eqref{ivhv1}}           \notag
\\
&
\begin{cases} 
& p,\,\, 2\le p\le\log d\,,
\\
&   \sqrt{(2p-1)+ (\log d-1)^2},\, \, \log d \le p\le \infty\,.
\end{cases}.
\end{align}
\end{thm}
This differs from inequality (3.2) of \cite{HT} and looks as it is better, but this is just because the assumption here is stronger: it  is in term of $P^2_p$ and in \cite{HT} the assumption of Lipschitzness is in terms of $K^2_\infty$.

\bigskip

But for $f\in M_{d\times d}$ we have obviously 
\begin{equation}
\label{pop}
\|\cdot\|_p \le d^{1/p} \|\cdot\|_{op}\,.
\end{equation}
Therefore inequality \eqref{m2} implies (again $\bE f=0$)
\begin{align}
\label{m3} 
&\!\!\!\!\!K^2_\infty \le 1 \Rightarrow \big(\bE\|f\|_{p}^p\big)^{1/p} \le^{\eqref{Kh}, \eqref{pop}} c_2\sqrt{p}\, d^{1/p} \,a(p)^{1/p}= c_2\sqrt{p} \,d^{1/p}e^{\beta(p)}  \le^{\eqref{ivhv1}}    \notag
\\
& 
\begin{cases} 
&c_2 d^{1/p} p^{3/2},\,\, 2\le p\le\log d\,,
\\
& c_2\,\sqrt{\log d} \, \sqrt{(2p-1)+ (\log d -1)^2},\, \, \log d \le p\le \infty\,.
\end{cases}
\end{align}
which is worse than (3.2) in \cite{HT}.

\bigskip


We have estimates  of $\sqrt{p}$-type for  large $p$  for all spaces $E$ of finite cotype simultaneously. So, this is not very surprising that they give worse estimates for concrete spaces of matrices.
Notice however, that using the Lipschitz condition in terms of gradient $P^2_p$ gives a new type of results \eqref{m2}.

Maybe it is worthwhile to mention that we consider all size $d$ matrix-functions on $\{-1,1\}^n$, while \cite{HT} is concerned with Hermitian matrix-functions.

Another difference is that we work only on Hamming cube, and \cite{HT} focuses on a wide class of Markov semigroups acting on matrices.

{}


\begin{thebibliography}{}






\bibitem[AW]{AW}{\sc Radosław Adamczak, Paweł Wolff}, 
{\em Concentration inequalities for non-Lipschitz functions with bounded derivatives of higher order}, 
Probab. Theory Related Fields 162 (2015), no. 3-4, 531--586.


\bibitem[CE]{CE}{\sc D. Cordero-Erasquin, A. Eskenazis}, {\em Discrete logarithmic Sobolev inequalities in Banach spaces}, J. London Math. Soc.,  
(2) 2024; 109:e12873., DOI: 10.1112/jlms.12873



\bibitem[GP]{GP} {\sc O. Gu\'edon and G. Paouris}, {\em  Concentration of mass on the Schatten classes}. Ann. Inst. H. Poincar\'e Probab. Statist., 43(1):87--99, 2007.

\bibitem[G]{G} {\sc L. Gross}, {\em Logarithmic Sobolev inequalities}. Amer. J. Math. 97 (1975), no. 4, 1061--1083. 

\bibitem[Ha]{Ha} {\sc U. Haagerup}, {\em The best constants in the Khintchine inequality}. Studia Math.,
70(3):231–283 (1982).

\bibitem[HT]{HT} {\sc De Huang, Joel A. Tropp}, {\em  Non-linear matrix concentration via semigroup methods}, Electronic J. Prob.v. 26 (2021), no.8, pp. 1--31.

\bibitem[HVNVW]{HVNVW}{\sc T. Hyt\"onen, J. Van Neerven, M. Veraar, L. Weiss}, Analysis in Banach spaces, vol. I, Martingales and Littlewood--Paley theory, Ergebnisse der Mathematik und ihrer Grenzgebiete, v. 63, Springer,
2016.

\bibitem[HVNVW2]{HVNVW2}{\sc T. Hyt\"onen, J. Van Neerven, M. Veraar, L. Weiss}, Analysis in Banach spaces, vol. II,  Ergebnisse der Mathematik und ihrer Grenzgebiete, v. 67, Springer,
2017.


\bibitem[IVHV]{IVHV}{\sc P. Ivanisvili, R. van Handel, and A. Volberg}, {\em Rademacher type and Enflo type coincide}, Ann. of Math. (2) 192
(2020), no. 2, 665--678.

\bibitem[L]{L} {\sc M. Ledoux}, {em The Concentration of Measure Phenomenon}, Volume 89 of Mathematical Surveys and
Monographs. American Mathematical Society, Providence (2001).

\bibitem[LO]{LO} {\sc M. Ledoux, K. Oleszkiewicz}, On measure concentration of vector-valued maps. Bull. Pol. Acad. Sci.
Math. 55(3), 261–278 (2007).

\bibitem[O]{O} {\sc K. Oleszkiewicz}, {\em Precise moment and tail bounds for Rademacher sums in
terms of weak parameters}. Israel J. Math., 203(1):429--443, 2014.

\bibitem[P]{P} {\sc G. Pisier}, {\em  Non commutative vector valued Lp spaces and completely $p$-summing
maps}. Ast\'erisque 247 (1998) SMF.



\bibitem[T]{T} {\sc M. Talagrand}, {\em Isoperimetry, logarithmic Sobolev inequality on the discrete cube and Margulis graph connectivity theorem}, 
GAFA, v.3, No. 3, 1993, pp. 295--314.

\bibitem[T1]{T1} {\sc  M. Talagrand}, {\em New concentration inequalities in product spaces}. Invent. Math. 126(3), 505--563 (1996).

\bibitem[TJ]{TJ} {\sc N. Tomczak-Jaegermann}, {\em  The moduli of smoothness and convexity and the Rademacher
averages of the trace classes $S_p$ ($1 \le p< \infty$).} Studia Math., 50:163–182, 1974.

\end{thebibliography}
\end{document}